\documentclass[12pt]{amsart}
\usepackage{amscd}
\advance\hoffset by -0.5 truecm
\usepackage{amssymb}
\newtheorem{Theorem}{Theorem}[section]
\newtheorem{Lemma}[Theorem]{Lemma}

\newtheorem{Remark}[Theorem]{Remark}

\def\Vol{\mbox{Vol}}
\def\V{\mbox{Var}}

\def\Z{{\mathbb Z}}
\def\R\re
\def\V{\bf V}

\def \re{{\mathbb R}}

\def \T{{\mathbb T}}

\def \0{\lambda_{0}}

\def \F{{\mathcal F}}

\begin{document}
\title[Collapsing manifolds]{Collapsing manifolds obtained by Kummer-type constructions}

\author[G. P. Paternain]{Gabriel P. Paternain}
\address{Department of Pure Mathematics and Mathematical Statistics,
University of Cambridge,
Cambridge CB3 0WB, England}
\email{g.p.paternain@dpmms.cam.ac.uk}

\author[J. Petean]{Jimmy Petean}
 \address{CIMAT  \\
          A.P. 402, 36000 \\
          Guanajuato. Gto. \\
          M\'exico.}
\email{jimmy@cimat.mx}

\thanks{J. Petean is supported by grant 37558-E of CONACYT}



\begin{abstract}
We construct $\F$-structures on a Bott manifold and on some other manifolds obtained by
Kummer-type constructions. We also prove that if $M=E\# X$ where $E$ is a fiber bundle with 
structure group $G$ and a fiber admitting a $G$-invariant metric of non-negative sectional 
curvature and $X$ admits an $\F$-structure with one trivial covering, then one can construct
on $M$ a sequence of metrics with sectional curvature uniformly bounded from below and 
volume tending to zero (i.e. $\Vol_K (M)=0$). As a corollary we prove that all the elements in the
Spin cobordism ring can be represented by manifolds $M$ with $\Vol_K (M)=0$.
  
\end{abstract}

\maketitle

\section{Introduction}

The present paper is concerned with the following invariants of closed smooth manifolds.
Given a closed connected smooth $n$-manifold $M$, let $\Vol(M,g)$ be the volume of a smooth
Riemannian metric $g$ and let $K_{g}$ be its sectional curvature.
Consider the following minimal volumes \cite{G2}:
\[\mbox{\rm MinVol}(M):=\inf_{g}\{\Vol(M,g):\;|K_{g}|\leq 1\}\]
and
\[\Vol_{K}(M):=\inf_{g}\{\Vol(M,g):\;K_{g}\geq -1\}.\]

A fundamental theorem of J. Cheeger and M. Gromov \cite{Cheeger} asserts that if $M$ admits a {\it polarized}
$\mathcal F$-{\it structure}, then $\mbox{\rm MinVol}(M)=0$, that is, we can collapse volume with
bounded sectional curvature. An $\F$-structure is a collection of tori acting on finite Galois coverings
of open subsets of the manifold. The actions are virtually effective, compatible with the finite group
of deck transformations of the coverings, and also compatible between themselves
on the overlap of the open subsets. These compatibility conditions ensure that $M$ is partitioned into orbits
which are flat manifolds.
The structure is said to be {\it polarized} if the dimension of the orbits is locally constant, in a certain
precise way. The $\mathcal F$-structure is said to be a $\mathcal T$-{\it structure} if the Galois coverings can all be taken
to be trivial. (We review these definitions in Section 2).

The vanishing of $\mbox{\rm MinVol}(M)$ implies, via Chern-Weil theory, that the 
Euler characteristic and all the Pontryagin numbers of $M$ are zero. The vanishing of $\Vol_{K}(M)$ implies that the {\it simplicial volume}
of $M$ is zero. Recall that the simplicial volume $||M||$ of a closed orientable manifold $M$ is defined as the
infimum of $\sum_{i}|r_{i}|$ where $r_i$ are the coefficient of a {\it real} cycle representing the fundamental
class of $M$. In fact, $||M||=0$ if there exists a sequence of metrics $g_i$ with
$\Vol(M,g_{i})\to 0$ and $\mbox{\rm Ric}_{g_{i}}\geq -(n-1)$ \cite{G2}. The vanishing of the simplicial volume
is the only known topological obstruction to $\Vol_{K}(M)=0$.  
The topology of closed orientable 3-manifolds with $\Vol_K =0$ has been determined
by T. Shioya and T. Yamaguchi \cite{SY1,SY2}. 
However, there could be smooth obstructions, as it occurs in dimension $4$. It is a striking consequence of 
Seiberg-Witten theory \cite{lebrun0,lebrun}  that if $M$ is a minimal compact complex surface of general type
then

\[\inf_{g}\left\{\Vol(M,g):\;\frac{s_{g}}{12}\geq -1\right\}=\frac{2\pi^2}{9}c^2_{1}(M), \]
where $s_{g}$ is the scalar curvature of $g$ and $c_{1}(M)$ is the first Chern class of $M$.

In \cite{PP0} we showed that if $M$ admits an arbitrary $\mathcal F$-structure, then $\Vol_{K}(M)=0$.
We also constructed $\mathcal T$-structures on several classes of manifolds, including compact complex elliptic
surfaces and any closed simply connected $5$-manifold. 

In \cite{Pe} the second author showed that if $M$ is a closed simply connected manifold of dimension
$\geq 5$, then

\[\inf_{g}\left\{\Vol(M,g):\;\frac{s_{g}}{n(n-1)}\geq -1\right\}=0.\]
More recently, C. Sung \cite{sung} extended this result by showing that
\[\inf_{g}\left\{\Vol(M,g):\;\lambda K_{g}+(1-\lambda)\frac{s_{g}}{n(n-1)}\geq -1\right\}=0\]
for any $\lambda\in [0,1)$.
All these results naturally raise the following question:

\medskip

\noindent {\bf Question 1.} {\it Let $M$ be a closed simply connected manifold of dimension
$n\geq 5$.
Is it true that $\mbox{\rm Vol}_{K}(M)=0$? }

\medskip

Since every closed simply connected 5-manifold admits a $\mathcal T$-structure, we know that Question 1
has a positive answer for $n=5$.
One can also speculate that if $M$ is odd-dimensional or if all its characteristic numbers are zero,
then $\mbox{\rm MinVol}(M)=0$. For a large class of highly connected manifolds, including certain exotic spheres 
which do not bound spin manifolds, this has been verified by C.Z. Tan \cite{tan}.

Here we show:

\medskip
\noindent {\bf Theorem A.} {\it Any closed connected spin manifold $M$
is spin cobordant to a manifold $N$ with $\mbox{\rm Vol}_{K}(N)=0$.
}

\medskip

We remark that Gromov and Lawson \cite{GL} proved that if $M$ and $N$ are spin cobordant
and $M$ is simply connected and of dimension $\geq 5$, then $M$ is obtained
from $N$ by performing surgery on spheres of codimension $\geq 3$.
If $M$ is not spin and in the same oriented cobordism class as $N$, then $M$ can also be obtained
from $N$ by performing surgery on spheres of codimension $\geq 3$ \cite[Proof of Theorem C]{GL}.
It follows easily from the description of the generators of $\Omega_{*}^{SO}$ as explained in
\cite[Proof of Theorem C]{GL} that every class in $\Omega_{*}^{SO}$ contains a manifold $N$
with $\mbox{\rm Vol}_{K}(N)=0$. Thus if one can show that the vanishing of $\Vol_{K}$ is invariant
under surgery on spheres of codimension $\geq 3$, Theorem A would imply that Question 1 has an
affirmative answer. We explicitly state this problem:

\medskip

\noindent {\bf Question 2.} {\it Let $M$ and $N$ be closed connected manifolds of dimension
$n\geq 5$. Suppose $\mbox{\rm Vol}_{K}(N)=0$ and $M$ is obtained
from $N$ by performing surgery on a sphere of codimension $\geq 3$. 
Is it true that $\mbox{\rm Vol}_{K}(M)=0$? }

\medskip

The proof of Theorem A is based on two novel ingredients. The first ingredient is the description
of an $\mathcal F$-structure on an 8-manifold $J_{8}$ with $\hat{A}$-genus equal to $1$ (a {\it Bott manifold}).
The manifold $J_{8}$ is one of the examples of 8-manifolds with special holonomy $\mbox{\rm Spin}(7)$
constructed by D. Joyce in \cite{J1}. In fact, we exhibit $\F$-structures on several manifolds with
special holonomy. These manifolds are all obtained by Kummer-type constructions, desingularizing a torus
orbifold. Besides the $K3$ and $J_{8}$ we exhibit polarized $\F$-structures on a Calabi-Yau 3-fold
with zero Euler characteristic and on a 7-manifold with holonomy $G_2$.
The case of the Calabi-Yau 3-fold is particularly interesting because it provides an example
of a closed simply connected 6-manifold manifold $M=X\# Y$ with $\mbox{\rm MinVol}(M)=0$, but with 
$\mbox{\rm MinVol}(X)$ and $\mbox{\rm MinVol}(Y)$ non-zero. Also $X$ is obtained from $S^6$ by surgery
on a complicated configuration of 3-spheres. We explain these observations in Remark \ref{cy}.

To our knowledge these examples of $\F$-structures on manifolds with special holonomy constitute
the first ones, besides the motivating case of flat manifolds, in which one can really appreciate the advantage
of the concept of $\F$-structure as opposed to the simpler concept of $\mathcal T$-structure.
On the other hand we do not know of an example of a manifold which carries an $\F$-structure but not
a $\mathcal T$-structure.

We observe that the existence of $\F$-structures on $K3$ and $J_{8}$ implies right away by our results in \cite{PP0}
that $\Vol_{K}(K3)=\Vol_{K}(J_{8})=0$. This fact seems to simplify some of the
proofs of the main results in \cite{Pe,sung}.

The second novel ingredient in the proof of Theorem A is the following:

\medskip

\noindent
{\bf Theorem B.} {\it Let $M=X \# E $ where $X$ is an $n$-manifold, $n\geq3$, which admits an $\F$-structure
with one trivial Galois covering and $E$ is the total space of a fiber bundle with fiber $F$ and structure
group $G$, where $F$ has a $G$-invariant metric of non-negative sectional curvature and $G$ is a compact connected
Lie group. Then $\mbox{\rm Vol}_K (M)=0$.} 

\medskip

\noindent
In fact, it seems plausible that a manifold like $E$ might
always admit an $\F$-structure. If that were the case, then we could deduce the vanishing of
$\Vol_{K}(X\# E)$ from \cite[Theorem 5.9]{PP0} which asserts that the connected sum of two manifolds with
$\F$-structures admits an $\F$-structure, provided that the $\F$-structures have at least one open set
with a trivial Galois covering. The proof of Theorem B uses a mixture of the collapsing techniques of
Cheeger and Gromov with computations of K. Fukaya  and T. Yamaguchi \cite[Theorem 0.18]{Yamaguchi} 
(see also \cite[Example 1.2]{Ya}) for the case
of bundles like $E$. One can see why Theorem B is important in the proof of Theorem A by recalling 
a celebrated result of S. Stolz \cite{S}: a closed spin manifold $M$ with zero $KO$-characteristic number
$\alpha(M)\in KO_{*}(\mbox{\rm point})$  is spin cobordant to the total space of a fibre bundle with fibre $\mathbb{HP}^{2}$
and structure group $\mbox{\rm PSp}(3)$. 

It is also important to note that in the condition $\Vol_K =0$
there is no implication on the diameter of the manifold; there are plenty of examples of manifolds which
are not almost non-negatively curved but nevertheless can be volume-collapsed with curvature bounded from
below (note for instance that for any closed manifold $M$, $\Vol_K (S^2 \times M)=0$). 
One does not expect spin manifolds with non-zero $\hat{A}$-genus like $K3$ and $J_8$ to be
almost non-negatively curved
(see \cite{lott}) nor do we expect such a conclusion for manifolds obtained as connected sums 
as in the previous paragraph.

Finally we mention that the existence of $\F$-structures also implies the vanishing of the {\it minimal entropy}
$\mbox{\rm h}(M)$ which is given by the infimum of the topological entropy of the geodesic flow of a metric $g$, as $g$ ranges over all metrics with volume $1$ \cite{PP0}. Most likely Theorem A holds also for $\mbox{\rm h}(M)$, but this would require to prove a result like Theorem B  for minimal entropy.
We do not pursue this issue here.
What we do obtain right away is $\mbox{\rm h}(K3)=\mbox{\rm h}(J_{8})=0$ as well as zero minimal
entropy for all the manifolds with special holonomy described in Section 3.

\medskip

{\it Acknowledgements:} The second author thanks the Department of Pure Mathematics and
Mathematical Statistics at the University of Cambridge and Trinity College for hospitality 
and financial support while this work was in progress.


\section{Preliminaries on $\F$-structures}

The notion of an $\F$-structure was first introduced by J. Cheeger and M. Gromov in 
\cite{G2,Cheeger, Cheeg}.
Although the definition is more elegantly expressed in terms of sheaves of tori and
local actions of them, in order to construct examples one usually uses equivalent definitions
in terms of an open cover of the manifold and torus actions on finite Galois coverings
of them. Since an important part of this article will be devoted to explicit constructions of
$\F$-structures we begin by giving the definition that we will use. 

An $\F$-structure on a closed manifold $M$ is given by the following data and conditions:

\begin{enumerate}
\item a finite open cover $\{U_{1},\dots,U_{N}\}$ of $M$;
\item $\pi_{i}:\widetilde{U}_{i}\to U_{i}$ is a finite Galois covering with group of deck transformations
$\Gamma_{i}$, $1\leq i\leq N$;
\item a smooth torus action with finite kernel of the $k_{i}$-dimensional torus, 
$\phi_{i}:T^{k_{i}}\to \mbox{\rm Diff}(\widetilde{U}_{i})$, $1\leq i\leq N$;
\item a homomorphism $\Psi_{i}:\Gamma_{i}\to \mbox{\rm Aut}(T^{k_{i}})$ such that
\[\gamma(\phi_{i}(t)(x))=\phi_{i}(\Psi_{i}(\gamma)(t))(\gamma x)\]
for all $\gamma\in \Gamma_{i}$, $t\in T^{k_{i}}$ and $x\in \widetilde{U}_{i}$;
\item for any finite sub-collection $\{U_{i_{1}},\dots,U_{i_{l}}\}$
such that $U_{i_{1}\dots i_{l}}:=U_{i_{1}}\cap\dots\cap U_{i_{l}}\neq\emptyset$
the following compatibility condition holds: let
$\widetilde{U}_{i_{1}\dots i_{l}}$ be the set of all points 
$(x_{i_{1}},\dots,x_{i_{l}})\in \widetilde{U}_{i_{1}}\times\dots\times\widetilde{U}_{i_{l}}$
such that $\pi_{i_{1}}(x_{i_{1}})=\dots=\pi_{i_{l}}(x_{i_{l}})$. The set
$\widetilde{U}_{i_{1}\dots i_{l}}$ covers $\pi_{i_{j}}^{-1}(U_{i_{1}\dots i_{l}})\subset \widetilde{U}_{i_{j}}$ for all
$1\leq j\leq l$. Then we require that $\phi_{i_{j}}$ leaves $\pi_{i_{j}}^{-1}(U_{i_{1}\dots i_{l}})$ invariant
and it lifts to an action on $\widetilde{U}_{i_{1}\dots i_{l}}$ such that all lifted actions commute.

\end{enumerate}

A subset  $S\subset M$ is called {\it invariant} if for any $y \in S\cap U_i $,
any  $x\in \widetilde{U_i}$ with $\pi_i (x)=y$ and any 
$t \in T^{k_i}$, we have $\pi_i (\phi_i (t) (x)) \in S$. 
The {\it orbit} of a point in $M$ is the minimal invariant set containing the point. 
The $\F$-structure is called {\it polarized} if given any subset $I\subset \{ 1,\dots ,N\}$,
if $U$ is the intersection of the $U_i's$ with $i\in I$ (assumed non-empty) and $V$ is the union of the $U_i's$
with $i\in I^c$, then  the dimension of the orbits is constant in $U-V$. The simplest case
in which the structure is polarized is when all the torus actions appearing in the definition are
locally free. An $\F$-structure is called a $\mathcal T$-structure if all the Galois coverings are
trivial.

We know of no previous example of an $\F$-structure which is not a 
$\mathcal T$-structure or the $\F$-structure on a flat manifold given by the holonomy covering. 
We will construct in the
following section a number of $\F$-structures on some interesting manifolds for which one of the
Galois coverings is not trivial. It is also not known to the authors if those manifolds admit any
$\mathcal T$-structure, except for the case of the $K3$ surface for which a $\mathcal T$-structure 
has been constructed in \cite[Theorem 5.10]{PP0}). It is important to note that in all 
these examples some of the Galois coverings (actually all but one) are trivial: 
this is needed in order to be able to construct
$\F$-structures on their connected sums as in \cite[Theorem 5.9]{PP0}.

\section{An exhibition of  $\F$-structures}

All the examples below will be desingularizations of torus orbifolds as in \cite{J0,J,J1}. 
The singularities
will be resolved using exclusively the {\it Eguchi-Hanson space} $X$ which is just a 4-manifold
diffeomorphic to $T^*S^2$ or to the blow up of ${\mathbb C}^2/\{\pm 1\}$. For the sake of clarity 
we will display the $\F$-structures in explicit cases although one can proceed in essentially
the same way in more general cases.

\subsection{The $K3$ surface}

The first exhibit in our gallery is the $K3$ surface.
It is already known that the $K3$ surface, as any compact complex elliptic surface, admits a 
$\mathcal T$-structure
\cite[Theorem 5.10]{PP0}. In this 
subsection we will put an $\F$-structure on the $K3$ surface described by the Kummer construction. This will give
a simpler way to see such a structure on the surface and will also provide a simple example for the constructions
in the following sections.

The $K3$ surface can be obtained as follows: first consider the map $J=-1 :\T^4 \rightarrow \T^4$, an involution of the
four-dimensional torus $\T^4=\re^4 /\Z^4$. The involution has 16 fixed points, those which have each coordinate
equal to 0 or 1/2. Let $P$ be the fixed point set of $J$. Then $\T^4 -P$ is invariant through $J$ and $J$ acts there 
without fixed points. The quotient $(\T^4 -P )/J$ is a smooth manifold with 16 ends diffeomorphic to 
${\mathbb R}{\mathbb P}^3 \times \re$. The $K3$ surface is obtained by attaching a copy of $X=T^* S^2$ 
to each one of the ends.

We describe an $\F$-structure on the $K3$ surface. Let $A_1$ be the circle action on $\T^4$ given by 
$x_{1}\mapsto x_{1}+\theta$ where $\theta\in \re/\Z$.
For $\varepsilon>0$ small, let 
\[\hat{U}:= \{ (x_1,\dots,x_4 )\in \T^4:\;d((x_2 ,x_3 ), (p_1 ,p_2 ))>\varepsilon /2\;
\mbox{\rm with}\; p_i =
0\;\mbox{\rm or}\;1/2 \}.\]
The open set $\hat{U}$ is invariant through $J$ and $A_1$. Let $U:=\hat{U}/J$. The set $U$ is an open 
subset of $K3$ and $\hat{U}$ is a $2:1$ covering of $U$. The homomorphism 
$\Psi:\{Id,J\}\to \mbox{\rm Aut}(S^1)$ given by
\[\Psi(J)(t)=-t\]
satisfies condition (4) in the definition of $\F$-structure.

Let 
\[V_1:=\{ (x_1,\dots,x_4 )\in \T^4:\; d((x_2 ,x_3 ),(0,0))<\varepsilon \}.\]
And define also $V_2$, $V_3$ and
$V_4$ replacing $(0,0)$ by $(0,1/2), (1/2,0)$ and $(1/2,1/2)$. Each $V_i$ is an open set of $\T^4$
diffeomorphic to $D^2 \times \T^2$. Let $B_i$ be the circle action on $V_i$ given by the obvious action on the
$D^2$-factor. The action $B_i$ commutes with the action $A_1$ and with $J$, and so it gives an action on
$(V_i -P)/J$, which we will also call $B_i$. 
The set $(V_i -P)/J$ is an open subset of $K3$ with 4 ends diffeomorphic to ${\mathbb R}{\mathbb P}^3$. 
On each of these ends the
action of $B_i$ is just the multiplication on two coordinates of $S^3$ which descends to an action on 
${\mathbb R}{\mathbb P}^3$.

Finally in a small metric ball around  each of the points in $P$ consider the usual Hopf action on $S^3$
which descends to ${\mathbb R}{\mathbb P}^3$. Call this action $C_i$; it commutes with $B_i$ and it extends to the whole $X$. 

The open subsets $U, V_1 ,V_2 ,V_3 ,V_4$ and the 16 copies of $X$ cover the whole $K3$ and the circle actions
$A_1 ,B_i, C_j$ define an $\F$-structure.

\subsection{An $\F$-structure on  $J_8$.} This would be the main attraction of our exhibition.
The existence of such a structure plays a key role in the proof of  Theorem A as we explained in
the introduction.

Consider the following involutions of the torus $\T^8=\re^8/\Z^8$:

\[\alpha(x_{1},\dots,x_{8})=(-x_{1},-x_{2},-x_{3},-x_{4},x_{5},x_{6},x_{7},x_{8}),\]
\[\beta(x_{1},\dots,x_{8})=\left(x_{1},x_{2},x_{3},x_{4},-x_{5},-x_{6},-x_{7},-x_{8}\right),\]
\[\gamma(x_{1},\dots,x_{8})=\left(\frac{1}{2}-x_{1},\frac{1}{2}-x_{2},x_{3},x_{4},\frac{1}{2}-x_{5},
\frac{1}{2}-x_{6},x_{7},x_{8}\right) \]
\[\delta(_{1},\dots,x_{8})=\left( -x_{1},x_{2},\frac{1}{2}-x_{3},x_{4},\frac{1}{2}-x_{5},x_{6},
\frac{1}{2}-x_{7},x_{8}\right) .\]

The fixed point sets of these involutions are:

\[S_{\alpha}=\{(p_{1},p_{2},p_{3},p_{4},x_{5},x_{6},x_{7},x_{8}):\;p_{i}=0\ \mbox{\rm or} \ 1/2\},\]
\[S_{\beta}=\{(x_{1},x_{2},x_{3},x_{4},p_{5},p_{6},p_{7},p_{8}):\;p_{i}=0\ \mbox{\rm or} \ 1/2 \},\]
\[S_{\gamma}=\{ (q_{1},q_{2},x_{3},x_{4},q_{5},q_{6},x_{7},x_{8})\;:q_{i}=1/4\ \mbox{\rm or}\ 3/4 \},\]
\[S_{\delta}=\{ (p_{1},x_{2},q_{3},x_{4},q_{5},x_{6},q_{7},x_{8}):\;p_{1}=0\ \mbox{\rm or} \ 1/2,\;q_{i}=1/4\
\mbox{\rm or} \ 3/4\}. \]

It is also easy to see that the involutions commute and generate a group $\Gamma$ isomorphic to
$(\Z_2 )^4$. The fixed point set of $\alpha \beta$ is $S_{\alpha} \cap S_{\beta}$ while all the other
elements of $\Gamma$ are fixed point free \cite[Example 1]{J1}. We will call $J_{8}$ the manifold obtained by resolving the
singularities of $\T^8 /\Gamma$ as below.

The set $S_{\gamma}$ is given by 16 copies of $\T^4$ and the group
generated by $\alpha, \beta$ and $\delta$ divides these into 2 groups of 8. Note also that $S_{\gamma}$
is disjoint from $S_{\alpha}$, $S_{\beta}$ and $S_{\delta}$. A neighbourhood of $S_{\gamma}$ then 
projects to $\T^8 /\Gamma$ onto 2 copies of $\T^4 \times (B^4 /{\pm 1})$. This singularity will be resolved
by replacing this open subset by $\T^4 \times X$. We will then have 2 open subsets $V_1 ,V_2$ on $J_{8}$ which
are diffeomorphic to $\T^4 \times X$. Note that if $Z$ is the zero section of $X=T^* S^2$ then 
$\T^4 \times (X -Z)$ lifts via the projection to 8 disjoint copies of $\T^4 \times (B^4 -\{ 0\} )$
in $\T^8$. The circle action on the last coordinate of the $\T^4$-factor will lift to each of the 
connected components as either as $\pm$ the circle action $A_{8}$ on the last coordinate of $\T^8$
($A_{8}$ is given by $x_{8}\mapsto x_{8}+\theta$ where $\theta\in \re/\Z$).

The set $S_{\delta}$ is resolved in the same way, producing open subsets $V_3$ and $V_4$ of $J_{8}$
where we also consider
the circle action induced by $A_{8}$.

The sets $S_{\alpha}$ and $S_{\beta}$ are also each of them equal to 16 copies of $\T^4$. The 16 $\T^4$'s corresponding to
$S_{\alpha}$ are divided into 4 groups of 4 by the action of the group generated by $\gamma$ and $\delta$; and
the same happens with $S_{\beta}$. The sets $S_{\alpha}$ and $S_{\beta}$ intersect in 256 points which are divided into
64 groups of 4 points by the action of $\gamma$ and $\delta$. Each of the components
of $S_{\alpha}$ is invariant by $\beta$ and $\beta$ acts on them as multiplication by $-1$; and the same is
true interchanging $\alpha$ and $\beta$. Therefore in resolving the singularities $S_{\alpha}$ we obtain
4 copies of $X\times K3$ and resolving $S_{\beta}$ we obtain 4 copies of $K3\times X$
(see the construction of $K3$ in the previous section). A connected component of the first type will
intersect a connected component of the second in 4 copies of $X \times X$.

Consider the open subset $\hat{U} \subset \T^8$ given by

\[\hat{U}= \{x\in \T^8 : d((x_{1},x_{2},x_{5},x_{6}),(q_{1},q_{2},q_{5},q_{6}))>\varepsilon /2, \]
\[ d((x_{1},x_{3},x_{5},x_{7}),(p_{1},q_{3},q_{5},q_{7}))>\varepsilon /2, \]
\[ d((x_{1},x_{3}),(p_{1},p_{3}))>\varepsilon /2 \ \mbox{\rm and}\]
\[d((x_{5},x_{7}),(p_{5},p_{7})>\varepsilon /2 , \ \mbox{\rm where}\]
\[ p_{i} =0 \ \mbox{\rm or} \ 1/2, \ q_{i}=1/4 \ \mbox{\rm or} \ 3/4 \} . \]

The set $\hat{U}$ is an
open subset of $\T^8$ which is disjoint from $S_{\alpha}$, $S_{\beta}$, $S_{\gamma}$ and $S_{\delta}$ and
is invariant through $A_8$ and $\Gamma$. Then $U=\hat{U}/\Gamma$ is an open subset of $J_{8}$ and
$\pi:\hat{U} \rightarrow U$  is an $8:1$ covering with deck transformation group $\Gamma$. We consider the
action $A_8$ on $\hat{U}$. The homomorphism $\Psi:\Gamma \rightarrow \mbox{\rm Aut}(S^1 )$ given by 
$\Psi(\alpha )=\Psi(\gamma )=\Psi(\delta )=Id$ and $\Psi(\beta)(t)=-t$ will satisfy the 
compatibility condition (4) in the definition of $\F$-structure. 

As we mentioned before, the action $A_8$ will give a circle action on the open subsets 
corresponding to the resolutions
of $S_{\gamma}$ and $S_{\delta}$; the preimage of the end of each $X \times \T^4$ 
(or $\T^4\times X$) under the projection
will have 8 connected components: in some of them we will have the action of $A_8$ and in some
of them we will have $-A_8$. Alternatively, one could also put the Hopf action on $X$.

Fix $p_{i}=0$ or $1/2$ and consider a small ball V around $(p_{1},p_{2},p_{3},p_{4})\in \T^4$. Then 
$V\times \T^4$ contains one of the components of $S_{\alpha}$. Note that  $\beta$ will act here as 
multiplication by $-1$ in the $\T^4$ factor. After resolving the singularities in the
projection $\T^8 /\Gamma$ we will have an
open subset of $J_{8}$  diffeomorphic to $X\times K3$. In $J_{8}$ there are 4 open subsets of this
type. Call them $W_{1},W_{2},W_{3}$ and $W_{4}$. And consider on them the Hopf action on the
$X$-factor. By considering the last 4 coordinates instead of the first, we obtain open
subsets $W_{5},W_{6},W_{7}$ and $W_{8}$ which contain 
(the resolutions of) the singularities $S_{\beta}$. And we put again the Hopf action on the $X$-factor.
For each $i\leq 4,j\geq 5$, $W_{i}$ will intersect $W_{j}$ in 4 copies of $X \times X$, which will
be invariant through both actions. These actions clearly commute.

Now consider

\[ {\hat{V}}_{1,3}= \{x\in \T^8 : \;d((x_{1},x_{3}),(p_{1},p_{3}))<\varepsilon,\;p_{i}=0 \ \mbox{\rm or} \ 1/2\} -S ,\]
where $S$ is the union of $S_{\alpha}$ and $S_{\beta}$. And let

\[{\hat{V}}_{5,7}=  \{x\in \T^8 : \;d((x_{5},x_{7}),(p_{5},p_{7}))<\varepsilon,\; \ p_{i}=0 \ \mbox{\rm or} \ 1/2\} -S .\]
Note that ${\hat{V}}_{1,3} $ and ${\hat{V}}_{5,7} $ are invariant through $\Gamma$.

Then ${\hat{V}}_{1,3} /\Gamma$ and ${\hat{V}}_{5,7} /\Gamma $ together with $U,V_{1},V_{2}, 
V_{3},V_{4},W_{1},\dots ,W_{8}$ cover
the whole of $J_{8}$. The set ${\hat{V}}_{1,3} /\Gamma$ is diffeomorphic 
to 
$$(D^2 \times \T^2 )/\{\pm 1\} \times \T^4 /\{\pm 1\}-S/\Gamma$$
and we put in this set the canonical circle action on the $D^2$-factor (it clearly leaves
$(D^2 \times \T^2 )/\{\pm 1\} \times \T^4 /\{\pm 1\}\cap S/\Gamma$ invariant).

This  action  lifts to the corresponding
open subset of $\hat{U}$ and commutes there with the $A_8$-action. 
We do the same thing for ${\hat{V}}_{5,7} /\Gamma $.
It is  easy to see that these
two actions will commute between themselves and with the actions on the $W_{i}$'s
in the corresponding intersections.

In summary: we covered $J_{8}$ with 15 open subsets. The sets $U$, ${\hat{V}}_{1,3} /\Gamma $ and${\hat{V}}_{5,7} /\Gamma $
are obtained as quotients of $\T^8$ by $\Gamma$ (away from fixed points). On $U$ we have the only non-trivial
covering for the $\F$-structure. The sets $V_{1},V_{2}, V_{3},V_{4}$ are open neighbourhoods of the resolutions of
the singularities $S_{\gamma}$ and $S_{\delta}$. They intersect $U$ but are disjoint from 
 ${\hat{V}}_{1,3} /\Gamma $ and ${\hat{V}}_{5,7} /\Gamma $. The sets $W_{1},\dots ,W_{8}$ cover the resolutions of the
singularities $S_{\alpha}$ and $S_{\beta}$. 
They only intersect ${\hat{V}}_{1,3} /\Gamma $ and ${\hat{V}}_{5,7} /\Gamma $.

\subsection{A polarized $\F$-structure on a closed 7-manifold with special holonomy $G_{2}$}
An interesting addition to the collection.

Let $(x_{1},\dots,x_{7})$ be coordinates in $\T=\re^7/\Z^7$ where $x_{i}\in \re/\Z$.
Let $\alpha,\beta$ and $\gamma$ be involutions of $\T^7$ defined by
\[\alpha(x_{1},\dots,x_{7})=(-x_{1},-x_{2},-x_{3},-x_{4},x_{5},x_{6},x_{7}),\]
\[\beta(x_{1},\dots,x_{7})=\left(-x_{1},\frac{1}{2}-x_{2},x_{3},x_{4},-x_{5},-x_{6},x_{7}\right),\]
\[\gamma(x_{1},\dots,x_{7})=\left(\frac{1}{2}-x_{1},x_{2},\frac{1}{2}-x_{3},x_{4},-x_{5},x_{6},-x_{7}\right).\]

One can easily check that $\alpha,\beta$ and $\gamma$ commute and hence they generate a group of isometries $\Gamma$ of the
flat torus $\T^7$ which is isomorphic to $\Z_{2}^3$. The following elementary properties of the action
of $\Gamma$ are proved in \cite{J0}.
The only non-trivial elements of $\Gamma$ that have fixed points are $\alpha, \beta$ and $\gamma$.
The fixed points of $\alpha$ in $\T^7$ are 16 copies of $\T^3$ and the group generated by $\beta$ and $\gamma$
acts freely on the set of 16 tori fixed by $\alpha$. Similarly the fixed points of $\beta$, $\gamma$ in $\T^7$
are each 16 copies of $\T^3$, and the groups $\langle \alpha,\gamma\rangle$ and $\langle \alpha,\beta\rangle$
act freely on the sets of 16 tori fixed by $\beta$ and $\gamma$ respectively.

Note that the 48 tori that make the set $S$ of points that are fixed by some non-trivial element in $\Gamma$ are all
disjoint. The singular set $S/\Gamma$ in the orbifold $\T^7/\Gamma$ is exactly the image of $S$ and 
consists of 12 copies of $\T^3$. For $\varepsilon$ sufficiently small each component of $S/\Gamma$ has a neighbourhood
isometric to $\T^3\times B^4_{\varepsilon}/\{\pm 1\}$. If we now resolve the singularities by replacing
 $\T^3\times B^4_{\varepsilon}/\{\pm 1\}$ by $\T^3\times X$, where $X$ is the Eguchi-Hanson space, we obtain
a closed simply connected 7-manifold $M$ which admits a family of metrics with holonomy $G_2$ \cite{J}.
The manifold $M$ has betti numbers $b_{2}=12$, $b_{3}=43$ and non-zero Pontryagin class $p_{1}(M)\in H^{4}(M,\Z)$.

Consider the following open sets in $\T^7$:

\[W^{\alpha}(\varepsilon):=\{(x_{1},\dots,x_{7})\in \T^7:\;d((x_{1},x_{2},x_{3}),(a_{1},a_{2},a_{3}))<\varepsilon\]
\[\mbox{\rm and}\;a_{i}=0,1/2,\;i=1,2,3\};\]
\[W^{\beta}(\varepsilon):=\{(x_{1},\dots,x_{7})\in \T^7:\;d((x_{1},x_{2},x_{5},x_{6}),(a_{1},a_{2},a_{3},a_{4}))<\varepsilon\]
\[\mbox{\rm and}\;a_{i}=0,1/2,\;i=1,5,6,\; a_{2}=1/4,3/4\};\]
\[W^{\gamma}(\varepsilon):=\{(x_{1},\dots,x_{7})\in \T^7:\;d((x_{1},x_{3},x_{5}),(a_{1},a_{3},a_{5}))<\varepsilon\]
\[\mbox{\rm and}\;a_{i}=1/4,3/4,\;i=1,3,\; a_{5}=0,1/2\}.\]

Note that for $\varepsilon$ small these sets are pairwise disjoint. Note also that $\Gamma$ leaves the 3 sets invariant.
The sets $W^{\alpha}(\varepsilon)$ and $W^{\delta}(\varepsilon)$ both have 8 connected components
which are copies of $\T^4\times B^{3}_{\varepsilon}$, while $W^{\beta}(\varepsilon)$ has 16 connected components
which are copies of $\T^3\times B^4_{\varepsilon}$. 

Now consider the open set $V$ in $\T^7$ given by the complement of the closure of 
$W^{\alpha}(\varepsilon/2)\cup W^{\beta}(\varepsilon/2)\cup W^{\beta}(\varepsilon/2)$. 
Observe that $\Gamma$ acts freely on $V$ and $V/\Gamma$ is an open set in $M$.

Let $\pi:\T^7\to \T^7/\Gamma$ be the projection map.
The set $\pi(W^{\alpha}(\varepsilon))$ will be an open set with 2 connected components which contain the 4 singular
3-tori corresponding to $\alpha$. When we resolve these singularities, $\pi(W^{\alpha}(\alpha))$ is modified
into an open set $U^{\alpha}(\varepsilon)$ of $M$. Similarly we obtain open sets $U^{\beta}(\varepsilon)$
and $U^{\delta}(\varepsilon)$ in $M$. By construction, the open sets $U^{\alpha}$, $U^{\beta}$, $U^{\delta}$ and $V/\Gamma$
cover $M$.

We now describe the torus actions. Let $A_{i}$ be the circle action on $\T^7$ given by $x_{i}\mapsto x_{i}+\theta$
where $\theta\in \re/\Z$. Observe that the action of $T^2$ given by $A_{4}\times A_{7}$ leaves invariant
$W^{\alpha}(\varepsilon)$, $W^{\beta}(\varepsilon)$, $W^{\delta}(\varepsilon)$ and $V$ (in fact, this property explains
why we have taken the above as relevant open sets).

Let us describe the circle action on $U^{\alpha}(\varepsilon)$. 
Note that $A_{7}$ commutes with $\alpha$ and $\beta$. Since $\delta$ swaps the connected components of 
$W^{\alpha}(\varepsilon)$ we can define using $A_{7}$ and $x_{7}\mapsto x_{7}-\theta$ a circle
action on $W^{\alpha}(\varepsilon)$ that will commute also with $\delta$.
This action descends to $\pi(W^{\alpha}(\varepsilon))$ and since the resolution of the singularities only affects
the first 4 coordinates we obtain a circle action $\phi_{\alpha}$ on $U^{\alpha}(\varepsilon)$.

Arguing similarly with $A_{4}$ and $\alpha$ for the sets $U^{\beta}(\varepsilon)$ and $U^{\delta}(\varepsilon)$ 
we obtain circle actions $\phi_{\beta}$ and $\phi_{\delta}$ on $U^{\beta}(\varepsilon)$ and $U^{\delta}(\varepsilon)$ 
respectively.

On the set $V/\Gamma$ we consider the action of $A_{4}\times A_{7}$ on $V$. The homomorphism 
$\Psi:\Gamma\to \mbox{\rm Aut}(T^2)$ given by
\[\Psi(\alpha)(t_{1},t_{2})=(-t_{1},t_{2}),\]
\[\Psi(\beta)(t_{1},t_{2})=(t_{1},t_{2}),\]
\[\Psi(\delta)(t_{1},t_{2})=(t_{1},-t_{2})\]
will clearly satisfy condition (4) in the definition of $\F$-structure.

Condition (5) in the definition of $\F$-structure follows quite easily from the fact that the only
possible overlaps are $(V/\Gamma)\cap U^{\alpha}$, $(V/\Gamma)\cap U^{\beta}$, $(V/\Gamma)\cap U^{\delta}$
and on them the actions lift and commute.

Finally, all actions are locally free and hence the $\F$-structure is polarized.

\subsection{A polarized $\F$-structure on a Calabi-Yau 3-fold}
\label{cy} The final exhibit.

Let $(x_{1},\dots,x_{6})$ be coordinates in $\T^6=\re^6/\Z^6$ where $x_{i}\in \re/\Z$.
Let $\alpha$ and $\beta$ be involutions of $\T^6$ defined by
\[\alpha(x_{1},\dots,x_{6})=(-x_{1},-x_{2},-x_{3},-x_{4},x_{5},x_{6}),\]
\[\beta(x_{1},\dots,x_{6})=\left(\frac{1}{2}-x_{1},-x_{2},x_{3},x_{4},-x_{5},-x_{6}\right).\]

As before, one can easily check that $\alpha$ and $\beta$ commute and hence they generate a group of isometries 
$\Gamma$ of the flat torus $\T^6$ which is isomorphic to $\Z_{2}^4$. 
The singular set of $\T^6/\Gamma$ consists of 16 copies of $\T^2$ each with a neighbourhood
isometric to $\T^2\times B^4_{\varepsilon}/\{\pm 1\}$. Desingularizing
these using the Eguchi-Hanson space yields a closed simply connected 6-manifold $M$
which carries a family of metrics with holonomy $SU(3)$ \cite[Example 2]{J}.
The manifold $M$ has betti numbers $b_{2}=19$, $b_{3}=40$ and thus zero Euler characteristic.

Consider the following open sets in $\T^6$:

\[W^{\alpha}(\varepsilon):=\{(x_{1},\dots,x_{6})\in \T^6:\;d((x_{1},x_{2},x_{3}),(a_{1},a_{2},a_{3}))<\varepsilon\]
\[\mbox{\rm and}\;a_{i}=0,1/2,\;i=1,2,3\};\]
\[W^{\beta}(\varepsilon):=\{(x_{1},\dots,x_{6})\in \T^6:\;d((x_{1},x_{2},x_{5}),(a_{1},a_{2},a_{5}))<\varepsilon\]
\[\mbox{\rm and}\;a_{i}=0,1/2,\;i=2,5,\; a_{1}=1/4,3/4\};\]

Note that for $\varepsilon$ small these sets are disjoint and $\Gamma$ leaves them invariant.
The sets $W^{\alpha}(\varepsilon)$ and $W^{\beta}(\varepsilon)$ both have 8 connected components
which are copies of $\T^3\times B^{3}_{\varepsilon}$.

Now consider the open set $V$ in $\T^6$ given by the complement of the closure of 
$W^{\alpha}(\varepsilon/2)\cup W^{\beta}(\varepsilon/2)$. 
Observe that $\Gamma$ acts freely on $V$ and $V/\Gamma$ is an open set in $M$.

Let $\pi:\T^6\to \T^6/\Gamma$ be the projection map.
The set $\pi(W^{\alpha}(\varepsilon))$ will be an open set with 4 connected components which contain the 8 singular
2-tori corresponding to $\alpha$. When we resolve these singularities, $\pi(W^{\alpha}(\alpha))$ is modified
into an open set $U^{\alpha}(\varepsilon)$ of $M$. Similarly we obtain an open set $U^{\beta}(\varepsilon)$
in $M$. By construction, the open sets $U^{\alpha}$, $U^{\beta}$ and $V/\Gamma$
cover $M$.

We now describe the torus actions. Observe that the action of $T^2$ given by $A_{4}\times A_{6}$ leaves invariant
$W^{\alpha}(\varepsilon)$, $W^{\beta}(\varepsilon)$ and $V$.

As before we get a circle action on $U^{\alpha}(\varepsilon)$ as follows. 
Note that $A_{6}$ commutes with $\alpha$. Since $\beta$ swaps the connected component of 
$W^{\alpha}(\varepsilon)$ we can define using $A_{6}$ and $x_{6}\mapsto x_{6}-\theta$ a circle
action on $W^{\alpha}(\varepsilon)$ that will commute also with $\beta$.
This action descends to $\pi(W^{\alpha}(\varepsilon))$ and since the resolution of the singularities only affects
the first 4 coordinates we obtain a circle action $\phi_{\alpha}$ on $U^{\alpha}(\varepsilon)$.

Arguing similarly with $A_{4}$ and $\alpha$ for the set $U^{\beta}(\varepsilon)$ 
we obtain a circle action $\phi_{\beta}$ on $U^{\beta}(\varepsilon)$.

On the set $V/\Gamma$ we consider the action of $A_{4}\times A_{6}$ on $V$. The homomorphism 
$\Psi:\Gamma\to \mbox{\rm Aut}(T^2)$ given by
\[\Psi(\alpha)(t_{1},t_{2})=(-t_{1},t_{2}),\]
\[\Psi(\beta)(t_{1},t_{2})=(t_{1},-t_{2}).\]
will clearly satisfy condition (4) in the definition of $\F$-structure.

Condition (5) in the definition of $\F$-structure follows quite easily from the fact that the only
possible overlaps are $(V/\Gamma)\cap U^{\alpha}$ and $(V/\Gamma)\cap U^{\beta}$
and on them the actions lift and commute.

Finally, all actions are locally free and hence the $\F$-structure is polarized.

\medskip

\begin{Remark}{\rm By Wall's splitting theorem for simply connected 6-manifolds \cite[Theorem 1]{W}, $M$ can be written
as $M=X\# Y$, where $X$ has $b_3=0$ and $Y$ is a connected sum of 20 copies of $S^3\times S^3$.
Clearly $X$ has positive Euler characteristic and $Y$ has negative Euler characteristic.
Hence the minimal volumes MinVol of $X$ and $Y$ are non-zero, but $\mbox{\rm MinVol}(M)=0$
since $M$ admits a polarized $\F$-structure. Note that $S^3\times S^3\#S^3\times S^3$ gives an
example of manifold with non-zero minimal volume, but $\mbox{\rm MinVol}(S^3\times S^3)=0$.
We conclude that the minimal volume does not behave well under connected sums.
According to Wall \cite[Theorem 2]{W}, the manifold $X$ can be obtained from $S^6$ by performing
surgery on a disjoint set of (framed) embedded 3-spheres. The spheres produce a link in $S^6$
and one can read off the cup form on $X$ (or $M$) from certain link invariants associated 
with the link (\cite[Theorem 4]{W}). The group $H^2(M,\re)$ admits a basis with 19 elements, 16 of which
come from the desingularization of $\T^2\times B^4_{\varepsilon}/\{\pm 1\}$, and the other 3 come
from $\Gamma$-invariant constant 
2-forms in $\T^6$ \cite[Section 2.3]{J}. If $e_i$ is any of these 19 elements, one can check that
$e_{i}^3=0$. Since $M$ is a Calabi-Yau 3-fold, the first Pontryagin class $p_1$ of $M$ is $-2c_{2}$
where $c_2$ is the second Chern class of $M$ and the latter must be non-zero. If we consider
$p_1$ as a linear form on $H^2(M,\re)$ we conclude that
$p_1(e_{i})$ is not zero for some $i$ and by \cite[Theorem 4]{W}, the class $e_i$
gives a rise to an embedding of $S^3$ in $S^6$ which is knotted. In principle one could compute
completely the cup form for this example as well as all the link invariants, but we 
do not pursue this matter here.}
\label{llave}
\end{Remark}

\section{Collapsing volume with curvature bounded from below}

In this section we prove Theorem B.

\begin{Lemma} Let $s:[2,4] \rightarrow \re$ be an non-decreasing smooth function which vanishes close
to 2, is equal to 1 close to 4, $s'\leq 1$ and $s'' \geq -2$. 
Let $h$ be a metric of non-negative sectional curvature
on a manifold $F$. For a small positive number $\delta$,
the metric $g=dt^2 + {\delta}^s h$ on $X= [2,4] \times F$ has sectional curvature bounded from
below by $-4  \log^2 (\delta )$. 
\end{Lemma}

\begin{proof} The metric $g$ is a warped product metric. Every plane $P\subset TX$ has a $g$- orthonormal
basis of the form $x+v,w$, where $v$ and $w$ are tangent to $F$ and $x$ is horizontal. 
Let $f={\delta}^s$. The sectional
curvature $K_g (P)$ is then computed by Bishop and O'Neill \cite[page 27]{Bishop}:

$$K_g (P) =-\frac{f''(t)}{f(t)} g(x,x) + \frac{K_h (Q) -(f'(t))^2}{f^2 (t)} g(v,v)$$
where $Q$ is the plane spanned by $v$ and $w$.
It follows that 

$$K_g (P) \geq -\frac{f''(t)}{f(t)}  -\frac{(f'(t))^2}{f^2 (t)}.$$

Now $f'=s' \log(\delta ) f$ and $f''= s''  \log(\delta ) f + s'^2 \log^2 (\delta ) f$. The lemma follows. 

\end{proof}

Let $G$ be a compact connected Lie group and $E$ be the total space of a fiber bundle with structure
group $G$ and fiber $F$ such that $F$ admits a $G$-invariant metric $g_F$ of non-negative sectional
curvature. Then there is a metric on $E$ for which the $G$-action is isometric and the fibers are
totally geodesic and have non-negative sectional curvature. If we do not consider the diameters of the 
manifolds, the computations of K. Fukaya  and T. Yamaguchi \cite[Theorem 0.18]{Yamaguchi} 
show that by shrinking the fibers one
obtains a sequence of metrics on $E$ with sectional curvature uniformly bounded from below and
collapsing volume. In \cite{PP0} the authors showed that the same is true for manifolds which admit
an $\F$-structure; the sequence of metrics is obtained by shrinking the orbits of the $\F$-structure
(away from the fixed points of the local torus actions). In the next lemma we will prove that
one can combine both cases at least when in one of the open subsets of the $\F$-structure the covering
is trivial (i.e. on one of the open subsets there is a torus acting). Note that one can assume
in this case that the torus acting is 1-dimensional.

\medskip

\noindent
{\bf Theorem B.} {\it Let $M =X \# E$ where $X$ is an n-manifold, $n\geq3$, which admits an $\F$-structure
with one trivial covering and $E$ is the total space of a fiber bundle with fiber  $F$ and structure
group $G$, where $F$ has a $G$-invariant metric of non-negative sectional curvature and $G$ is a 
compact connected Lie Group.
Then $\mbox{\rm Vol}_K (M)=0$.}

\medskip

\begin{proof} Let ${\bf g}$ be a metric on $E$  such that all the fibers are totally geodesic
and isometric to some fixed Riemannian manifold 
$(F,g_F )$ of non-negative sectional curvature (on which
$G$ acts by isometries). 
We will moreover assume
that over some disk $B^k$ in the base space, the bundle is trivial and the metric is a 
nice Riemannian
product; namely, there is an open subset of $E$ which is diffeomorphic to $B^k \times F$ and
the restriction of ${\bf g}$ is the product of $g_F$  and a  metric
of non-negative sectional curvature on $B^k$ for which the end is isometric to
$[1,5]\times S^{k-1}$. 

Inside $B^k$ consider the sets $V_1 =[4,5]\times S^{k-1}$, $V_2=[2,4]\times S^{k-1}$,
$V_3 =[1,2]\times S^{k-1}$ and $B_2 =B(0,2)\subset B^k$.

To perform the connected sum of $X$ and $E$ we will pick a point $e\in V_3 \times F$ and
a point $x$ in the open subset of $X$ with a trivial covering, $x$ a regular point for the
corresponding circle action. 
Pick any linear circle action on $S^{k-1}$ and use it to define a circle action on $B^k \times F$. 
Assume that $e$ lies on a regular orbit for this action. On the connected sum of $X$ with $B_2 \times F$
we construct an $\F$-structure as in \cite[Theorem 5.9]{PP0}: consider tubular neighbourhoods of the orbits
through $x$ and $e$,
$S^1 \times B_X$ and $S^1 \times B_E$. One then uses that $S^1 \times B^{k-1} \#
S^1 \times B^{k-1}$ is diffeomorphic to $S^1 \times B^{k-1} - S^{k-2} \times B^2$ to 
construct an $\F$-structure  on the connected sum of $S^1 \times B_X$ and 
$S^1 \times B_E$ which will match the original actions on both connected components of
the boundary.

Define a Riemannian metric ${\bf h}$ on $M$ so that it coincides with 
${\bf g}$ away from $V_3 \times F$ and is invariant 
for the $\F$-structure above. This is achieved by the same procedure as in \cite[Lemma 1.3]{Cheeger},
averaging by the local torus actions a metric which coincides with ${\bf g}$ away from $V_3 \times F$.

\vspace{.5cm}

Fix $\delta$ small. We are going to construct a  metric ${\bf h}_{\delta}$ on $M$. 

\vspace{.5cm}

Let us denote by 
$X^*$ the part of $M$ where we have the $\F$-structure; namely, the connected sum of 
$B_2 \times F$ with $X$. On $X^*$ we proceed as in \cite{PP0}: we 
consider $X^* \times T^N$ (for some appropriate $N$) and define a diagonal polarized $\F$-structure
on  $X^* \times T^N$.  Namely, if for the $\F$-structure we have open subsets
$U_1 ,\dots ,U_l$ where the corresponding local torus action is not locally free we put 
$N=k_1 +\dots +k_l$ and on each $U_i \times T^N$ we consider the corresponding diagonal action.
On $X^* \times T^N$ all the local actions are locally free and we can perform the construction 
of Cheeger and Gromov \cite[Section 3]{Cheeger}
obtaining a sequence of metrics with bounded sectional curvature. 
Recall that the sequence of metrics is obtained by first multiplying the metric by
$\log^2 (\delta )$ and then multiplying by $\delta$ (or more precisely by some appropriate function
of $\delta$) in the direction of the orbits.
Then we take the quotient by the
$T^N$-action. The result is a sequence of metrics ${\bf h}_{\delta}^0$ on $X^*$ so that 
$\Vol(X^* ,{\bf h}_{\delta}^0 ) \rightarrow  0$ while the sectional curvature remains bounded from 
below (this computation is carried out in \cite{PP0}). Moreover on the boundary the metric
${\bf h}_{\delta}^0$ will be the Riemannian quotient by the $S^1$ action of 

$$\log^2 (\delta )
\left( g_F \times dt^2 |_{[2-\varepsilon ,2]}\times q_{\delta} \right) ,$$ 

\noindent
where $q_{\delta}$ is the metric on $S^{k-1} \times S^1$  obtained from the product metric by
multiplying by $\delta$ in the tangent space of the diagonal action.

Away from $B^k \times F$ we put the metric ${\bf h}_{\delta}^1$ obtained from 
$\log^2 (\delta ) {\bf g}$ by multiplying the metric
by $\delta$ in the directions of the fiber. The fact that the curvature remains bounded from
below follows from the O'Neill formulas 
as shown in \cite[Theorem 0.18]{Yamaguchi} (see also \cite[Example 1.2]{Ya}). 
It is also clear that the volume
of this region also goes to 0 with $\delta$.

Finally we have to join  ${\bf h}_{\delta}^0$ and ${\bf h}_{\delta}^1$ along $[2,5]\times S^{k-1} \times F$. 
On $V_1 \times F$ we leave the metric equal to $\delta \log^{2} (\delta )  g_F$ on the
$F$-factor and we modify the metric on $[4,5]\times S^{k-1}$: 
let $s_1 :[4,5]\rightarrow \re$ be a smooth function which is equal to 0 near 5 and equal to 1
near 4,  
consider the product metric on $[4,5]\times  S^{k-1} \times
S^1$ and multiply by ${\delta}^{s_1}$ the tangent space of the diagonal action. 
Multiplying by $\log^2 (\delta )$ and taking the quotient
by the $S^1$-factor we obtain a metric which glues well with ${\bf h}_{\delta}^1$ on one extreme and
on the other extreme will coincide with ${\bf h}_{\delta}^0$ on the $V_1$-factor. 
It is trivial that the volume
of this region will go to 0 with $\delta$. Cheeger and Gromov \cite[Theorem 3.1]{Cheeger}
prove that the sequence of metrics
on $V_1 \times S^1$ have bounded sectional curvature as $\delta$ tends to 0. Since the sectional
curvature does not decrease by Riemannian submersions, the sectional curvature of the metric
on $V_1$ will remain bounded from below. 

Now we modify the metric on $V_2 \times F$: here we prolong the metric defined at
$t=4$ in the previous paragraph on the $S^{k-1}$-factor
and on the $[2,4]\times F$-factor we put the metric
$ \log^2 (\delta ) (dt^2 + {\delta}^{s_2 } g_F )$, where $s_2 : [2,4]\rightarrow \re$ is a smooth function as in
the Lemma 4.1.
The fact that in this region the curvature remains bounded 
from below follows from  the previous lemma. The fact that the volume collapses is clear
since the volume of the $S^{k-1}$-factor is already collapsed.
This metric will glue well with the previous one where $s_2 =1$ and will glue well with 
${\bf h}_{\delta}^0$ when $s_2 =0$.

\end{proof}

\section{Proof of  Theorem A}

We need to show that every element of the Spin cobordism ring $\Omega_*^{\mbox{\rm Spin}}$ can be 
represented
by a manifold $N$ with $\Vol_K (N)=0$. There is a surjective ring homomorphism 
$\alpha :\Omega_*^{\mbox{\rm Spin}} \rightarrow KO_* (\mbox{\rm point})$. S. Stolz \cite{S} proved that every 
element in the kernel of $\alpha$ is represented by the total space of a fiber bundle with
fiber $\mathbb{HP}^{2}$
and structure group $\mbox{\rm PSp}(3)$. 

On the other hand if $a\in KO_n$ then there exists a closed 
spin manifold $X$ with an $\F$-structure with one trivial Galois covering 
such that $\alpha [X] =a$: if $B$ is a Bott
manifold, multiplication by $\alpha ([B])$ gives an isomorphism between
$KO_n$ and $KO_{n+8}$. Then the construction in Subsection 3.2 proves the claim in dimension 
$\geq 8$ and the low dimensions
can be easily dealt with by hand as in \cite[Theorem 2]{Pe}. 

Therefore every element of $\Omega_*^{\mbox{\rm Spin}}$
can be represented by a manifold of the form $N=E\# X$ where $X$ admits an $\F$-structure with
one trivial Galois covering and $E$ is the total space of a fiber bundle with
fiber $\mathbb{HP}^{2}$
and structure group $\mbox{\rm PSp}(3)$. The theorem follows from Theorem B.

\end{document}